\documentclass{amsproc}
\usepackage{amssymb}
\usepackage{amsmath}
\usepackage{amsfonts}

\theoremstyle{plain}

\numberwithin{equation}{section}
\input{tcilatex}

\begin{document}
\title[New Proofs on Properties of an Orthogonal Decomposition]{New Proofs
on Properties of an Orthogonal Decomposition of a Hilbert Space}
\author{Dejenie A. Lakew}
\address{John Tyler Community College\\
Department of Mathematics}
\email{dlakew@jtcc.edu}
\urladdr{http://www.jtcc.edu}
\date{October 21, 2015}

\begin{abstract}
We establish new and different kinds of proofs of properties that arise due
to the orthogonal decomposition of the Hilbert space, including projections,
over the unit interval of one dimension. We also see angles between
functions, particularly between those which are non zero constant mulitples
of each other and between functions from $A^{2}\left( \Omega \right) $ and $%
D\left( W_{0}^{1,2}\left( \Omega \right) \right) $.
\end{abstract}

\maketitle

$\overline{\ \ \ \ \ \ \ \ \ \ \ \ \ \ \ \ \ \ \ \ \ \ \ \ \ \ \ \ \ \ \ \ \ 
}$

$\ \ \ $

\textbf{Notations}.\ \ \ 

\ 

Let $\Omega =\left[ 0,1\right] $, $\partial \Omega =\{0,1\}$ \ \ \ \ 

\ \ \ 

$\oplus $\ :\ \ \ \ \ Set orthogonal direct sum

\ \ 

$\uplus $ \ :\ \ \ \ \ Elements direct sum from mutually orthogonal sets \ 

\ \ \ \ \ \ \ \ \ \ \ \ \ \ \ \ \ \ \ \ \ \ \ \ \ \ \ \ \ 

$D:=\frac{d}{dx}$

\ \ \ \ \ \ \ \ \ \ \ \ \ 

As customarily we define

\ \ \ \ 

$(i)$ \ \ The Hilbert space of square integrable functions over $\Omega $ \ 
\begin{equation*}
\tciLaplace ^{2}\left( \Omega \right) =\{f:\Omega \longrightarrow 
\mathbb{R}
,\text{ measurable and \ }\dint\limits_{\Omega }f^{2}dx<\infty \}
\end{equation*}

\ \ \ \ 

$(ii)$ \ The Sobolev space 
\begin{equation*}
W^{1,2}\left( \Omega \right) =\{f\in \tciLaplace ^{2}\left( \Omega \right)
:Df\in \tciLaplace ^{2}\left( \Omega \right) \}
\end{equation*}%
and

\ \ 

$(iii)$ \ the traceless Sobolev space 
\begin{equation*}
W_{0}^{1,2}\left( \Omega \right) =\{f\in W^{1,2}\left( \Omega \right) :\tau
^{+}f=f_{\mid \partial \Omega }=0\}
\end{equation*}

where $\tau ^{+}$ is trace to the boundary from the inside of the interval
of a function in a particular function space. The Hilbet space $\tciLaplace
^{2}\left( \Omega \right) $ is an inner product space with inner product 
\begin{equation*}
\langle ,\rangle :\tciLaplace ^{2}\left( \Omega \right) \times \tciLaplace
^{2}\left( \Omega \right) \longrightarrow 
\mathbb{R}%
\end{equation*}%
defined by 
\begin{equation*}
\langle f,g\rangle =\dint\limits_{\Omega }f\left( x\right) g\left( x\right)
dx
\end{equation*}

\ \ \ \ \ \ \ \ \ \ 

Thus norm of a function in the Hilbert space is defined by

\ \ \ 
\begin{equation*}
\parallel f\parallel _{\tciLaplace ^{2}\left( \Omega \right) }:=\sqrt{%
\left\langle f,f\right\rangle }=\left( \int_{\Omega }f^{2}dx\right) ^{\frac{1%
}{2}}
\end{equation*}

\ \ \ 

The theme of this short article is to present new ways of looking at
partitioning a function from a Hilbert space as a direct sum of a function
that is well behaving and constant over the domain and an other function
from the derivative of the Sobolev space. The problem which looks apparent
and to be addressed is, how every function from the Hilbert space can be
decomposed in the way described, as a sum of a function that is
differentiable with zero derivative over the entire domain and some other
summand from a $D$-image of a Sobolev space. This is because not every
function in the Hilbert space is well behaving or smooth enough so that we
differentiate with no difficulty. Some of them are discontinuous, or blow up
at interior or boundary points and what the new proof addresses is these
pitfalls by creating well functioning and behaving partitioning of such
functions so that we differentiate with the derivatives to be smooth and
zero. We construct the constant function out of such a function by taking
its integral value over the domain which is always a finite number by the
fact that a function in the Hilbert space is always Lebesgue integrable from
the argument,

\ \ 
\begin{eqnarray*}
&\parallel &f\parallel _{\tciLaplace ^{1}\left( \Omega \right)
}=\int_{\Omega }\mid f\mid dx \\
&\leq &\left( \int_{\Omega }\mid f\mid ^{2}dx\right) ^{\frac{1}{2}}\left(
\int_{\Omega }dx\right) ^{\frac{1}{2}} \\
&=&\parallel f\parallel _{\tciLaplace ^{2}\left( \Omega \right) }<\infty
\end{eqnarray*}%
The next process is to construct the second direct summoned from $D\left(
W_{0}^{1,2}\left( \Omega \right) \right) $.

\ \ \ \ 

With respect to the defined inner product we have the following orthogonal
decomposition of the space given by

\ \ \ \ 

\textbf{Proposition 1}. \ (\textit{Orthogonal Decomposition) The Hilbert
space }$\tciLaplace ^{2}\left( \Omega \right) $ has an orthogonal
decomposition given by

\ \ \ \ \ \ 
\begin{equation*}
\ \tciLaplace ^{2}\left( \Omega \right) =A^{2}\left( \Omega \right) \oplus
D\left( W_{0}^{1,2}\left( \Omega \right) \right)
\end{equation*}

\ \ 

\textbf{Proof}. \ \ 

\ 

$\bigskip $We need to show two things

$(i)$ \ $A^{2}\left( \Omega \right) \cap D\left( W_{0}^{1,2}\left( \Omega
\right) \right) =\{0\}$

\ \ 

$(ii)$ $\ $Every function $f$ in the Hilbert space $\tciLaplace ^{2}\left(
\Omega \right) $ has a unique representation as a direct sum of functions
from the two subspaces in a unique way, i.e., 
\begin{equation*}
\forall f\in \tciLaplace ^{2}\left( \Omega \right) ,\exists g\in A^{2}\left(
\Omega \right) \text{ \ \ and \ \ \ }\exists h\in D\left( W_{0}^{1,2}\left(
\Omega \right) \right)
\end{equation*}

such that%
\begin{equation*}
f=g\uplus h.
\end{equation*}

\ \ 

The proof of $\left( i\right) $ is available in $\left[ 1\right] $

\ \ \ \ 

And the proof of $\left( ii\right) $ follows from the following proposition.

\ \ \ 

\textbf{Proposition 2. }Let $f\in \tciLaplace ^{2}\left( \Omega \right) $,
then 
\begin{equation*}
\exists h\in W_{0}^{1,2}\left( \Omega \right)
\end{equation*}%
such that%
\begin{equation*}
f(x)=\int_{\Omega }f(x)dx\uplus Dh(x)
\end{equation*}

\ \ 

\textbf{Proof.} Consider 
\begin{equation*}
h\left( x\right) =\int_{0}^{x}f(t)dt-x\int_{\Omega }f(x)dx
\end{equation*}%
First let us show that $h\in W_{0}^{1,2}\left( \Omega \right) $. Clearly 
\begin{equation*}
\tau ^{+}h=0\text{ on }\partial \Omega \text{\ \ as }h\left( 0\right)
=0=h\left( 1\right)
\end{equation*}%
and $h\in \tciLaplace ^{2}\left( \Omega \right) $ for the fact that 
\begin{equation*}
\int_{0}^{x}f(t)dt\text{ \ \ \ and \ \ }x\int_{\Omega }f(x)dx
\end{equation*}%
\ are both in $\tciLaplace ^{2}\left( \Omega \right) $

\ 

Also 
\begin{equation*}
Dh(x)=f(x)-\int_{\Omega }f(x)dx\in \tciLaplace ^{2}\left( \Omega \right)
\end{equation*}

\begin{equation*}
\therefore \ \ \ \ h\in W_{0}^{1,2}\left( \Omega \right)
\end{equation*}

But%
\begin{equation*}
\int_{\Omega }f(x)dx\uplus \underset{\underset{Dh}{\parallel }}{\underbrace{%
\left( f(x)-\int_{\Omega }f(x)dx\right) }}=f(x)
\end{equation*}

Hence follows the proposition.

\ 

This is a very important result as it indicates how we construct the
constant function out of the smooth or how ill behaving a function may be
that remains square integrable. Now $(ii)$ of \textit{proposition 1} easily
follows.

\ \ 

$(ii)$ Let $f\in \tciLaplace ^{2}\left( \Omega \right) $. Then following 
\textit{proposition 2}, we have made all the necessary ingredients for the
splitting process

\ \ \ 
\begin{equation*}
g=\int_{\Omega }f(x)dx\ \ \ \ \text{and \ \ \ \ \ \ }\eta =Dh(x)
\end{equation*}%
where $h$ is given above, so that%
\begin{equation*}
f=g\uplus \eta
\end{equation*}%
\ with 
\begin{equation*}
g\in A^{2}\left( \Omega \right) \text{,\ \ }g\text{ is a constant and hence }%
Dg=0,\text{ }g\in \tciLaplace ^{2}\left( \Omega \right)
\end{equation*}%
and 
\begin{equation*}
\eta \in D\left( W_{0}^{1,2}\left( \Omega \right) \right)
\end{equation*}%
which proves $(ii)$.

\ 

\textbf{Definition}. \ Due to the orthogonal decomposition, there are two
orthogonal projections

\ 
\begin{equation*}
P:\tciLaplace ^{2}\left( \Omega \right) \longrightarrow A^{2}\left( \Omega
\right)
\end{equation*}%
\ and 
\begin{equation*}
Q:\tciLaplace ^{2}\left( \Omega \right) \longrightarrow D\left(
W_{0}^{1,2}\left( \Omega \right) \right)
\end{equation*}%
with 
\begin{equation*}
Q=I-P
\end{equation*}%
where $I$ is the identity operator.

\ 

\textbf{Proposition 3}. \ \ $\ \ \forall f\in \tciLaplace ^{2}\left( \Omega
\right) $ \ \ we have \ \ 
\begin{equation*}
\left\langle P\left( f\right) ,Q\left( f\right) \right\rangle =0
\end{equation*}%
\ 

\textbf{Proof }. From the decomposition result of proposition 2 
\begin{equation*}
P(f)=\int_{\Omega }f(x)dx
\end{equation*}%
\ and 
\begin{equation*}
Q(f)=f(x)-\int_{\Omega }f(x)dx
\end{equation*}%
we have 
\begin{eqnarray*}
\left\langle P(f),Q(f)\right\rangle &=&\int_{\Omega }P(f)Q(f)dx \\
&=&\int_{\Omega }\int_{\Omega }f(t)dt\left( f(x)-\int_{\Omega }f(t)dt\right)
dx \\
&=&\int_{\Omega }\left( \int_{\Omega }f(t)dt\right) f(x)dx-\int_{\Omega
}\left( \int_{\Omega }f(t)dt\right) ^{2}dx \\
&=&0\text{, since }\int_{\Omega }dx=1
\end{eqnarray*}

\ \ 

\textbf{Proposition 4. }\ \ \ We have the following properties

\ \ \ \ \ \ \ \ \ \ \ \ \ \ \ \ \ \ \ \ \ \ \ \ \ \ \ \ \ 

$(i)$ \ \ \ $P\circ Q=Q\circ P=0$

\ \ \ \ \ \ \ \ \ \ \ \ \ \ \ \ \ \ \ 

$(ii)$ \ \ $P^{2}=P$

\ \ \ \ \ \ \ \ \ \ \ \ \ \ \ \ \ \ \ \ \ \ 

$(iii)$ $\ Q^{2}=Q$

\ \ \ \ \ \ \ \ \ \ \ \ \ \ \ \ \ \ \ \ \ \ \ \ \ 

That is $P$ and $Q$ are \textit{idempotent} operators.

\ \ \ \ \ \ \ \ \ \ \ \ \ \ \ \ \ \ \ \ \ \ \ \ \ \ \ 

\textbf{Proof. }

\ 

$\left( i\right) $ Let $f\in \tciLaplace ^{2}\left( \Omega \right) $

\ 

then from \textbf{proposition 2 }%
\begin{equation*}
P(f)=\int_{\Omega }f(x)dx
\end{equation*}%
\ and 
\begin{equation*}
\ Q(f)=f(x)-\int_{\Omega }f(x)dx
\end{equation*}%
and then 
\begin{eqnarray*}
P(Q(f)) &=&P\left( \left( f(x)-\int_{\Omega }f(x)dx\right) \right) \\
&=&\int_{\Omega }\left( f(x)-\int_{\Omega }f(y)dy\right) dx \\
&& \\
&=&0\text{ \ }
\end{eqnarray*}

\ \ and%
\begin{eqnarray*}
Q(P(f)) &=&Q\left( \int_{\Omega }f(x)dx\right) \\
&=&\int_{\Omega }f(x)dx-\int_{\Omega }\left( \int_{\Omega }f(y)dy\right) dx
\\
&& \\
&=&0
\end{eqnarray*}

Hence 
\begin{equation*}
Q(P)=P(Q)=0
\end{equation*}

\bigskip

$\left( ii\right) $ follows from%
\begin{equation*}
f=P(f)\uplus Q(f)
\end{equation*}
and then applying $P$ on $f$ again we have 
\begin{eqnarray*}
P(f) &=&P\left( P(f)\uplus Q(f)\right) \\
&=&P^{2}(f)\uplus P(Q(f)) \\
&=&P^{2}(f)
\end{eqnarray*}%
as 
\begin{equation*}
P(Q(f))=0
\end{equation*}

That is%
\begin{equation*}
P^{2}=P
\end{equation*}

$\left( iii\right) $ Similarly let $f\in \tciLaplace ^{2}\left( \Omega
\right) $. Then%
\begin{equation*}
f=P\left( f\right) \uplus Q\left( f\right)
\end{equation*}%
and therefore 
\begin{equation*}
Q\left( f\right) =Q\left( P\left( f\right) \uplus Q\left( f\right) \right)
=QP\left( f\right) \uplus Q^{2}\left( f\right)
\end{equation*}

But 
\begin{equation*}
Q\left( P\left( f\right) \right) =0
\end{equation*}%
and hence 
\begin{equation*}
Q^{2}=Q
\end{equation*}

$\ \ \ \ \ \ \ \ \ \ \ \ \ \ \ \ \ \ \ \ \ \ \ \ \ \ \ \ \ \ \ \ $

\textbf{Proposition 6}. $\ \ \forall f\in \tciLaplace ^{2}\left( \Omega
\right) $%
\begin{equation*}
\dint\limits_{\Omega }Q\left( f\right) dx=0
\end{equation*}

\textbf{Proof}. For $f\in \tciLaplace ^{2}\left( \Omega \right) $ using the
decomposition \textit{proposition 2},%
\begin{equation*}
Q\left( f\right) =f(x)-\int_{\Omega }f(t)dt
\end{equation*}%
we have 
\begin{eqnarray*}
\int_{\Omega }Q(f)(x)dx &=&\int_{\Omega }\left( f(x)-\int_{\Omega
}f(t)dt\right) dx \\
&=&\int_{\Omega }f(x)dx-\int_{\Omega }\int_{\Omega }f(t)dtdx \\
&=&0
\end{eqnarray*}%
as 
\begin{equation*}
\int_{\Omega }dx=1
\end{equation*}

\ \ \ 

\textbf{Proposition 7}. Let $f\in \tciLaplace ^{2}\left( \Omega \right) \cap
C^{1}\left( \Omega \right) $. Then 
\begin{equation*}
\exists x_{0}\in \Omega ^{\text{int}}:Q\left( f\left( x_{0}\right) \right) =0
\end{equation*}

\ \ 

\textbf{Proof}. Suppose $\nexists x_{0}\in \Omega ^{\text{int}}$: $Q\left(
f\left( x_{0}\right) \right) \neq 0$. Then either 
\begin{equation*}
Q\left( f\left( x\right) \right) >0\text{ \ \ on \ \ }\Omega ^{\text{int}%
}\Longrightarrow \dint\limits_{\Omega }Q\left( f\right) dx>0
\end{equation*}

or 
\begin{equation*}
Q\left( f\left( x\right) \right) <0\text{ \ on \ \ }\Omega ^{\text{int}%
}\Longrightarrow \dint\limits_{\Omega }Q\left( f\right) dx<0
\end{equation*}

Thus%
\begin{equation*}
\dint\limits_{\Omega }Q\left( f\right) dx\neq 0
\end{equation*}%
which is a contradiction to \textit{proposition 6}.

\ \ \ \ 

Hence the proposition follows.

\ \ \ \ \ 

\textbf{Corollary 8}. For $f\in \tciLaplace ^{2}\left( \Omega \right) $, we
have 
\begin{equation*}
\int_{\Omega }P\left( f\right) dx=\int_{\Omega }fdx
\end{equation*}

\ \ \ \ \ \ 

\textbf{Proof}. The result follows from \textbf{proposition 6}.

\ \ \ \ \ \ \ 

\ \ \ \ \ \ 

\textbf{Proposition 9}. For $f\in \tciLaplace ^{2}\left( \Omega \right) $,
if 
\begin{equation*}
P(f)=0,\text{ then }\int_{\Omega }Q(f)(x)dx=0
\end{equation*}

\ 

\textbf{Proof}. \ 
\begin{equation*}
P(f)=\int_{\Omega }f(x)dx=0
\end{equation*}%
Then 
\begin{equation*}
Q(f)=f(x)-P(f)=f(x)
\end{equation*}%
and hence 
\begin{eqnarray*}
\int_{\Omega }Q(f(x))dx &=&\int_{\Omega }f(x)dx \\
&=&P(f) \\
&=&0
\end{eqnarray*}

\ \ 

\textbf{Proposition 10}. 
\begin{equation*}
P,Q:\tciLaplace ^{2}\left( \Omega \right) \longrightarrow \tciLaplace
^{2}\left( \Omega \right)
\end{equation*}%
with 
\begin{equation*}
\ker P=\{f\in \tciLaplace ^{2}\left( \Omega \right) :\int_{\Omega }f(x)dx=0\}
\end{equation*}%
and 
\begin{equation*}
\ker Q=\{f\in \tciLaplace ^{2}\left( \Omega \right) :f(x)=\int_{\Omega
}f(x)dx\}
\end{equation*}

\textbf{Proof}. From the decomposition \textit{proposition 2}, 
\begin{equation*}
P(f)=\int_{\Omega }f(x)dx
\end{equation*}%
and hence 
\begin{eqnarray*}
P(f) &=&\int_{\Omega }f(x)dx=0 \\
&\Longrightarrow &\text{ \ }\int_{\Omega }f(x)dx=0
\end{eqnarray*}

Also%
\begin{eqnarray*}
Q(f) &=&f(x)-\int_{\Omega }f(x)dx=0 \\
&\Longrightarrow &\text{ \ \ }f(x)=\int_{\Omega }f(x)dx
\end{eqnarray*}

\ \ \ \ \ \ 

$\ \ \ \ \ \ \ \ \ \ \ \ \ \ \ \ \ \ \ \ \ \ \ \ \ \ \ \ \ \ \ \ \ \ \ \ \ \
\ \ \ \ \ \ \ \ \ \ \ \ $

$\overline{\ \ \ \ \ \ \ \ \ \ \ \ \ \ \ \ \ \ \ \ \ \ \ \ \ \ \ \ \ \ \ \ \
\ \ \ \ }$

\ \ 

We will look at few illustrations on these results.

\ \ \ \ 

\textbf{Example 1. \ \ \ }\ \ For \ $f\left( x\right) =x,$ 
\begin{equation*}
P(x)=\frac{1}{2}\text{ \ \ and \ \ }Q(x)=x-\frac{1}{2}
\end{equation*}
so that 
\begin{equation*}
x=\frac{1}{2}\uplus \left( x-\frac{1}{2}\right)
\end{equation*}%
$\ \ $

\textbf{Proof}. \ \ \ 

\ \ \ \ \ 

Indeed from \textit{proposition 2},

\begin{equation*}
\ P(f)=\int_{\Omega }xdx=\frac{1}{2}\ \ \ \text{and \ \ }Q(f)=x-\frac{1}{2}
\end{equation*}

In a similar construction we have

\ 

\textbf{Example 2.}\ \ For $f\left( x\right) =x^{2}$, $\ $%
\begin{equation*}
P\left( f\right) =\int_{\Omega }x^{2}dx=\frac{1}{3}\ \ \ \ \text{and}\ \
Q\left( f\right) =x^{2}-\frac{1}{3}
\end{equation*}%
so that 
\begin{equation*}
x^{2}=\frac{1}{3}\uplus \left( x^{2}-\frac{1}{3}\right)
\end{equation*}

\ 

With similar procedures we decompose the follwing functions,

\ \ \ \ \ \ \ 

\textbf{Example 3}. $\ f(x)=\mid x-\frac{1}{2}\mid $, 
\begin{equation*}
P(f(x))=\frac{1}{4}\text{ \ \ and \ }Q(f(x))=\mid x-\frac{1}{2}\mid -\frac{1%
}{4}
\end{equation*}
so that 
\begin{equation*}
\mid x-\frac{1}{2}\mid =\frac{1}{4}\uplus \left( \mid x-\frac{1}{2}\mid -%
\frac{1}{4}\right)
\end{equation*}%
\ \ \ 

\textbf{Example 4}. For $\gamma \in 
\mathbb{R}
\setminus \{0\}$ and for $f(x)=e^{\gamma x}$, 
\begin{equation*}
P(e^{\gamma x})=\frac{1}{\gamma }\left( e^{\gamma }-1\right) ,\text{ \ }%
Q(e^{\gamma x})=e^{\gamma x}+\frac{1}{\gamma }-\frac{e^{\gamma }}{\gamma }
\end{equation*}
so that

$\ \ $%
\begin{equation*}
\ e^{\gamma x}=\frac{1}{\gamma }\left( e^{\gamma }-1\right) \uplus \left(
e^{\gamma x}+\frac{1}{\gamma }-\frac{e^{\gamma }}{\gamma }\right)
\end{equation*}

\textbf{Example 5}.$\ \ f(x)=\cos \gamma x,$ 
\begin{equation*}
P(\cos \gamma x)=\frac{\sin \gamma }{\gamma },\text{ \ }\ Q\left( \cos
\gamma x\right) =\cos \gamma x-\frac{\sin \gamma }{\gamma }
\end{equation*}

\ 

so that 
\begin{equation*}
\cos \gamma x=\frac{\sin \gamma }{\gamma }\uplus \left( \cos \gamma x-\frac{%
\sin \gamma }{\gamma }\right)
\end{equation*}

for $\gamma \in 
\mathbb{R}
\setminus \{0\}$

\ \ 

The case for functions that are not nice enough, this procedure works very
well.

\ \ 

\textbf{Example 6}. Let $f(x)=\left( x-\frac{1}{2}\right) ^{\frac{2}{3}}$,
clearly $f$ is not differentiable over $\Omega $ but decomposable as%
\begin{eqnarray*}
P(f) &=&\int_{\Omega }\left( x-\frac{1}{2}\right) ^{\frac{2}{3}}dx \\
&=&\frac{3}{5\sqrt[3]{4}}\in A^{2}\left( \Omega \right)
\end{eqnarray*}%
and 
\begin{equation*}
Q(f)=\left( x-\frac{1}{2}\right) ^{\frac{2}{3}}-\frac{3}{5\sqrt[3]{4}}\in
D\left( W_{0}^{1,2}\left( \Omega \right) \right)
\end{equation*}%
so that 
\begin{equation*}
f(x)=\frac{3}{5\sqrt[3]{4}}\uplus \left( \left( x-\frac{1}{2}\right) ^{\frac{%
2}{3}}-\frac{3}{5\sqrt[3]{4}}\right)
\end{equation*}

\ \ 

\textbf{Example 7}. Let $f\in H^{2}(\Omega )$ which is discontinuous at
several points and badly behaving. Then $f\in \tciLaplace ^{1}\left( \Omega
\right) $ and thus let 
\begin{equation*}
P(f)=\int_{\Omega }f(x)dx
\end{equation*}%
is finite and hence in $A^{2}\left( \Omega \right) $ and 
\begin{equation*}
Q(f)=f(x)-P(f)\in D\left( W_{0}^{1,2}\left( \Omega \right) \right)
\end{equation*}%
so that%
\begin{equation*}
f(x)=P(f)\uplus Q(f)
\end{equation*}

\ \ \ \ \ \ 

$\overline{\ \ \ \ \ \ \ \ \ \ \ \ \ \ \ \ \ \ \ \ \ \ \ \ \ \ \ \ \ \ \ \ \
\ \ \ \ \ \ \ }$\ \ \ 

\ \ \ \ \ 

We can also discuss about \textit{angles} between functions in $\tciLaplace
^{2}\left( \Omega \right) $ as it is an inner product space.

\ \ \ \ \ 

\textbf{Definition}. For $f,g\in \tciLaplace ^{2}\left( \Omega \right) $, we
define the angle $\theta $ between them as 
\begin{equation*}
\theta =\cos ^{-1}\left( \frac{\left\langle f,g\right\rangle }{\parallel
f\parallel \parallel g\parallel }\right)
\end{equation*}

\ \ 

Clearly $\forall f\in A^{2}\left( \Omega \right) ,\forall g\in D\left(
W_{0}^{1,2}\left( \Omega \right) \right) $ with both non zero, the angle
between them is%
\begin{equation*}
\theta =\frac{\pi }{2}\text{ \ \ as \ \ }\left\langle f,g\right\rangle =0
\end{equation*}

\ \ 

\textbf{Proposition 11}. For two non zero functions $f$ and $g$ in $%
\tciLaplace ^{2}\left( \Omega \right) $ with $g=\lambda f$ where $\lambda $
is a non zero constant, the angle $\theta $ between $f$ and $g$ is 
\begin{equation*}
\theta =\theta \langle (f,g)=\left\{ 
\begin{array}{c}
0\text{, for }\lambda >0 \\ 
\pi \text{, for }\lambda <0%
\end{array}%
\right.
\end{equation*}

\textbf{Proof.} 
\begin{eqnarray*}
\theta &=&\theta \langle (f,g) \\
&=&\cos ^{-1}\left( \frac{\left\langle f,\lambda f\right\rangle }{\mid
\lambda \mid \parallel f\parallel ^{2}}\right) \\
&=&\cos ^{-1}\left( \frac{\lambda \parallel f\parallel ^{2}}{\mid \lambda
\mid \parallel f\parallel ^{2}}\right) \\
&=&\cos ^{-1}\left( \frac{\lambda }{\mid \lambda \mid }\right) \\
&=&\left\{ 
\begin{array}{c}
0\text{, for }\lambda >0 \\ 
\pi \text{, for }\lambda <0%
\end{array}%
\right.
\end{eqnarray*}

\textbf{Example 6}. Compute the angle between $f(x)=x$ and $g(x)=x^{2}$

Solution. 
\begin{equation*}
\parallel f\parallel =\left( \int_{\Omega }x^{2}dx\right) ^{\frac{1}{2}}=%
\frac{1}{\sqrt{3}}
\end{equation*}
\ and 
\begin{equation*}
\parallel g\parallel =\left( \int_{\Omega }x^{4}dx\right) ^{\frac{1}{2}}=%
\frac{1}{\sqrt{5}}
\end{equation*}
Also 
\begin{equation*}
\left\langle f,g\right\rangle =\int_{\Omega }x^{3}dx=\frac{1}{4}
\end{equation*}

Thus the angle between $f$ and $g$ is 
\begin{equation*}
\theta =\cos ^{-1}\left( \frac{1/4}{1/\sqrt{15}}\right) =\cos ^{-1}\left( 
\frac{\sqrt{15}}{4}\right)
\end{equation*}

\textbf{Exercises}. Find the angle between the following pair of functions

\ 

$(a)$ $f(x)=\sin x$ and $g(x)=\cos x$

\ 

$(b)$ $f(x)=e-1$ and $g(x)=e^{x}+1-e$

\ 

$(c)$ $f(x)=e^{x}$ and $g(x)=e^{-x}$

\ \ \ \ \ 

\ \ 

$\overline{\ \ \ \ \ \ \ \ \ \ \ \ \ \ \ \ \ \ \ \ \ \ \ \ \ \ \ \ \ \ \ \ \
\ \ \ \ \ \ }$

\ \ \ \ \ \ \ \ \ 

\textbf{Future} \textbf{research works}.

\ \ \ \ \ 

Are the following decompositions valid ?

\ 

$\left( i\right) $ $W^{1,2}\left( \Omega \right) =A^{1,2}\left( \Omega
\right) \oplus D^{2}(W_{0}^{3,2}\left( \Omega \right) $

\ 

$\left( ii\right) $ $W^{k-1,2}\left( \Omega \right) =A^{k,2}\left( \Omega
\right) \oplus D^{k}\left( W_{0}^{2k-1,2}\left( \Omega \right) \right) $

\ 

with $D^{2}:=\frac{d^{2}}{dx^{2}}$ \ and $D^{k}:=\frac{d^{k}}{dx^{k}}$ \ and 
$k\in 
\mathbb{N}
$

\ 

$\left( iii\right) $ Conjecture: For $f\in \tciLaplace ^{2}\left( \Omega
\right) \cap C^{1}\left( \Omega \right) $, $Qf$ \ has some kind of symmetry
over $\Omega $ either in terms of function property or in terms of integral
values.

\ \ 

For instance for $f(x)=e^{x},$%
\begin{equation*}
\exists x_{0}=\ln \left( e-1\right) \in \left( 0,1\right) =\Omega ^{\text{int%
}}:Q\left( f(x_{0}))\right) =0
\end{equation*}%
so that 
\begin{equation*}
\int_{0}^{\ln \left( e-1\right) }Q(f(x))dx=-\int_{\ln \left( e-1\right)
}^{1}Q(f(x))dx
\end{equation*}

the very reason why $\ $%
\begin{equation*}
\int_{\Omega }Q(f(x))dx=0
\end{equation*}

\ \ \ \ 

\textbf{N.B}. When I was presenting a seminar on this topic in the
Department of Mathematics and Applied Mathematics regular analysis, logic
and physics seminar, at Virginia Commonwealth University, faculty members
from the audience indicated to me that this probably might be an
intermediate value theorem version in this setting and I give credit to them
for pointing that out.

\bigskip

\ \ \ 

\ \ \ \ \ \ \ \ \ \ \ \ \ \ \ \ \ \ \ \ \ \ \ \ \ \ \ \ \

\ 


\begin{thebibliography}{9}
\bibitem{dej1} Dejenie A. Lakew, On Orthogonal Decomposition of a Hilbert
Space $\tciLaplace ^{2}\left( \Omega \right) $, International Journal of
Mathematics and Computer Science, No$.1,10(2015)27-37$

\bibitem{dej2} Dejenie A. Lakew, Mulugeta Alemayehu, Clifford Analysis on
Orlicz-Sobolev Spaces, $arXiv:1409.8380v1$

\bibitem{dejray} Dejenie A. Lakew, John Ryan, Complete Function Systems and
Decomposition Results Arising in Clifford Analysis, Computational Methods
and Function Theory, No. $1(2002)215-228$
\end{thebibliography}
\end{document}